\documentstyle[12pt]{article}
\input {cyracc.def}
\font\tencyr=wncyr7 scaled \magstep1
\def\rus{\tencyr\cyracc}

\newcommand{\reb}[1]{\mbox{\bf  (\ref{#1})}}
\newcommand{\re}[1]{\mbox{\rm  (\ref{#1})}}

\newenvironment{proof}
{\noindent {\sl Proof.}\quad }{\hfill
$\square$ \vskip1.1ex\noindent }

\newenvironment{proof*}
{\noindent {\sl Proof.}\quad }{\hfill
$\square$}

\renewcommand{\theequation}{\thesection .\arabic{equation}}

\catcode`\@=\active
\catcode`\@=11
\def\@eqnnum{\hbox to .01pt{}\rlap{\bf \hskip -\displaywidth(\theequation)}}
\catcode`\@=12

\newenvironment{s}[1]
{ \vskip1.2ex \refstepcounter{equation}
\noindent {\bf \theequation\quad #1.} \begin{sl}}{\end{sl}
\vskip1.1ex\noindent }

\newenvironment{rem}[1]
{ \vskip1.2ex \refstepcounter{equation}
\noindent {\bf \theequation\quad {#1}.} }{ \vskip1.1ex\noindent }

\newenvironment{subs}[1]
{\vskip1.2ex \refstepcounter{equation}
\noindent {\bf (\theequation)\quad #1.} }{\quad}

\newcommand {\sekt}[1]
{{\vskip2.5ex\refstepcounter{section}\setcounter{equation}{0}
\noindent\large \bf \thesection\quad #1
\vskip1.5ex\noindent}}

\newcommand {\ah}{{\frak a}}

\newcommand {\ce}{{\frak c}}

\newcommand {\g}{{\frak g}}
\newcommand {\h}{{\frak h}}
\newcommand {\ka}{{\frak k}}
\newcommand {\el}{{\frak l}}
\newcommand {\me}{{\frak m}}
\newcommand {\n}{{\frak n}}

\newcommand {\p}{{\frak p}}
\newcommand {\es}{{\frak s}}
\newcommand {\te}{{\frak t}}
\newcommand {\tri}{{\frak sl}_2}

\newcommand {\z}{{\frak z}}
\newcommand {\esi}{\varepsilon}
\newcommand {\ap}{\alpha}

\newcommand {\lb}{\lambda}

\newcommand {\tih}{\tilde h}

\newcommand {\V}{{\Bbb V}}

\newcommand {\hb}{{\bf h}}
\newcommand {\fb}{{\bf f}}

\newcommand {\eb}{{\bf{e}}}

\newcommand {\psumma}{{\in}\!\!\!{+}}
\newcommand {\mathrm}{\rm\mbox}
\newcommand {\isom}{\stackrel{\sim}{\longrightarrow}}
\newcommand {\ad}{{\mathrm{ad\,}}}
\newcommand {\Ad}{{\mathrm{Ad\,}}}

\newcommand {\cox}{{\mathrm{Cxt}}}

\newcommand {\hot}{{\mathrm{ht\,}}}

\newcommand {\Ker}{{\mathrm{Ker\,}}}
\newcommand {\Ima}{{\mathrm{Im\,}}}
\newcommand {\rk}{{\mathrm{rk\,}}}

\newcommand {\GR}[2]{{\mathrm{{\bf #1}}}_{#2}}

\newcommand {\vno}[1]{\vskip#1 ex\noindent}
\newcommand {\rar}{\rightarrow}

\newcommand {\beq}{\begin{equation}}
\newcommand {\eeq}{\end{equation}}

\newfam\Bbbfam\newfam\eufam%
\font\Bbbfont=msbm10 scaled 1200%
\font\olala=msam10 scaled 1200%
\font\frak=eufm10 scaled 1400%
\font\Bbbsmallfont=msbm7 scaled 1000%
\textfont\Bbbfam=\Bbbfont\scriptfont\Bbbfam=\Bbbsmallfont
\font\euzw=eufm10 scaled 1200%
\font\euac=eufm7 scaled 1200%
\textfont\eufam=\euzw\scriptfont\eufam=\euac%
\def\frak{\fam\eufam}%
\def\Bbb{\fam\Bbbfam}%
\def\bbb{\fam\Bbbfam\scriptstyle}

\def\square{\hbox {\olala\char"03}}
\def\Bbbk{\hbox {\Bbbfont\char'174}}

\def\ltimes{\hbox {\Bbbfont\char'156}}

\begin{document}
\setlength{\parskip}{2pt plus 4pt minus 0pt}
\hfill {\scriptsize June 7, 1999} 
\vskip1ex

\noindent
{\Large \bf Nilpotent pairs in semisimple Lie algebras 
 \vno{1.3}and their characteristics}
\bigskip \\
{\bf Dmitri I. Panyushev}

\medskip
\smallskip

\noindent{\large \bf Introduction}
\vno{2}%
In a recent article \cite{vitya}, V.\,Ginzburg introduced and studied in
depth the notion of a {\sl principal nilpotent pair\/} in a semisimple
Lie algebra $\g$. He also obtained several results for more general pairs.
As a next step, we considered in \cite{pCM} {\sl almost principal
nilpotent pairs\/}.
The aim of this paper is to make a contribution to the general theory of
nilpotent pairs. Roughly speaking, a nilpotent pair $\eb=(e_1,e_2)$ consists
of two commuting elements in $\g$ that can independently be contracted to
the origin (see precise definition in sect.\,1).
A principal nilpotent pair is a double counterpart of a regular (=principal)
nilpotent element. Consequently, the theory of nilpotent pairs should stand out
as double counterpart of the theory of nilpotent orbits. As the cornerstone
of the latter is the Morozov--Jacobson theorem
and the concept of a characteristic,
the primary goal is to
realize to which extent these can be generalized to the double setting.
\\[.7ex]
The fundamental distinction of the double situation is that there is no
general analogue of the Morozov--Jacobson theorem.
What remains true is that any nilpotent pair has a characteristic
$\hb=(h_1,h_2)$, which is unique within to conjugacy \reb{charact}.
Hence characteristics can be used to study further properties of nilpotent
pairs. Generalizing Dynkin's approach in \cite{EBD}, we prove that the number
of $G$-orbits of characteristics of nilpotent pairs is finite \reb{finita}
and provide some information about the numerical labels $\ap_j(h_i)$, where
$\{\ap_1,\dots,\ap_n\}$ is a suitably chosen set of simple roots of $\g$
\reb{labels}. Since the number of $G$-orbits of
nilpotent pairs is infinite (see \cite[5.5]{vitya}), one encounters a
challenging problem to restore somehow a ``one-to-one" correspondence between
nilpotent pairs and characteristics. Our solution is that
we introduce in sect.\,2 {\sl wonderful\/} (nilpotent) pairs. The definition
is given in terms of the bi-grading of $\g$ determined by $\hb$ and involves
only $\hb$-eigenspaces with integral eigenvalues \reb{wonder}.
We prove that if
$\eb,\eb'$ are two wonderful pairs with the same characteristic
$\hb$, then $\eb$ and $\eb'$ are $Z_G(\hb)$-conjugate \reb{wond3}.
This implies that there are finitely many $G$-orbits of wonderful pairs.
On the other hand, specific classes of nilpotent pairs
considered in \cite{vitya},\cite{wir},\cite{pCM} are wonderful.
\\[.7ex]
In section 3, properties of several classes of wonderful pairs are studied.
For instance, we call a nilpotent pair $\eb$ {\sl even\/}, if $\dim\z_\g(\hb)
=\dim\z_\g(\eb)$. It is a natural analogue of an even nilpotent element.
We show that if $\eb$ is even, then it is wonderful and the eigenvalues of
$\hb$ in $\g$ are integral \reb{even}.
\\[.7ex]
In section 4, we describe characteristics of principal and almost principal
nilpotent pairs \reb{pr-char} and indicate the relationship between
eigenvalues of $\hb$ in $\z_\g(\eb)$ and the exponents of the
corresponding Levi subalgebras $\z_\g(h_i)$ \reb{exponents}.
\\[.7ex]
What is still lacking is a true understanding of possible fractional eigenvalues
of $\hb$ for arbitrary wonderful pairs. Unlike the ordinary theory
(=\,theory of nilpotent orbits),
these eigenvalues may have very large denominators \reb{denom}.
\\[.7ex]
The ground field $\Bbbk$ is algebraically closed and of
characteristic zero. Algebraic groups are denoted by capital Roman letters,
while their Lie algebras by the corresponding small Gothic letters.
Throughout, 
$\g$ is a semisimple Lie algebra and $G$ is its adjoint group.
For any set $M\subset\g$, we let $\z_\g(M)$ (resp. $Z_G(M)$) denote
the centralizer of $M$ in $\g$ (resp. in $G$) and $M^\perp$ the orthogonal
complement in $\g$ with respect to the Killing form.
For $M=\{a,\dots,z\}$, we simply write $\z_g(a,\dots,z)$ or $Z_G(a,\dots,z)$.
If $N\subset G$, then $Z_G(N)$ stands for the centralizer of $N$ in $G$.
For $x\in\g$ and $s\in G$, we write $s{\cdot}x$
in place of $(\Ad s)x$. Write $A^o$ for the identity component of an algebraic
group $A$. \par
$\langle a,\dots,z\rangle$ is the linear span of
elements of a vector space; \par
${\Bbb P}=\{0,1,2,\dots\}$, ${\Bbb N}=\{1,2,\dots\}$.
\\
Our general reference for algebraic groups is Vinberg--Onishchik's book
\cite{vion}. For nilpotent orbits consult \cite{SS}, \cite[ch.\,6]{t41},
and \cite{CoMc}.
\vno{.8}%
{\small {\bf Acknowledgements.}
This research was supported in part by RFFI Grant {\rus N0}\, 98--01--00598.
I would like to thank T.\,Levasseur and P.\,Torasso for arranging my visit
to Poitiers (France), where a portion of this work was done.}

\sekt{Characteristics and their properties\nopagebreak}%
Let us begin with a definition, which is due to V.\,Ginzburg.
\begin{rem}{Definition} \label{nil-p}
A pair $\eb=(e_1,e_2) \in \g\times\g$
is said to be {\it nilpotent\/} in $\g$, if \\
(i) $[e_1, e_2]=0$  and
(ii) for any $(t_1,t_2)\in \Bbbk^*\times \Bbbk^*$,  there
exists $g=g(t_1,t_2)\in G$  such that
$(\,t_1e_1,\,t_2e_2\,)=(g{\cdot}e_1,\,g{\cdot}e_2\,).$
\end{rem}%
Obviously, then both $e_1$ and $e_2$ are nilpotent elements of $\g$.
Note however that a nilpotent pair is not the same as a commuting
pair of nilpotent elements. A nilpotent pair is said to be {\it trivial\/},
if one of $e_i$'s is equal to 0.
It was shown in \cite{vitya} (see also \ref{charact}(1) below)
that condition (ii) is equivalent to
the following one: there exists
a pair of semisimple elements $\hb=(h_1, h_2)\in\g\times\g$ such that
$\ad h_1$ and $\ad h_2$ have rational eigenvalues and
\beq [h_1, h_2]=0,\quad [h_i, e_j]=\delta_{ij} e_j \quad
(i,j\in \{1,2\})\ .  \label{comrel}
\eeq
Ginzburg called such a pair an associated semisimple pair. He also proved that
an associated semisimple is unique up to conjugacy, if $\eb$ is a
``pre-distinguished" nilpotent
pair \cite[sect.\,5]{vitya}. We shall prove that, after a slight modification
of the definition, the result actually holds for all nilpotent pairs. For
this reason, we prefer to use the classical term
`characteristic' introduced by E.\,B.\,Dynkin in \cite{EBD}.
\begin{rem}{Definition}
A pair of semisimple elements $\hb\in\g\times\g$ is called a
{\it characteristic\/} of $\eb$, if
it satisfies commutation relations \re{comrel} and $\{h_1,h_2\}
\subset\z_\g(\eb)^\perp$.
\end{rem}%
Recall that if $A$ is
a linear algebraic group, then there exists an algebraic Levi decomposition
$A=A^{red}\ltimes A^{nil}$, where $A^{nil}$
is the unipotent radical and $A^{red}$ is a reductive Levi
subgroup of $A$ (see \cite[ch.\,6]{vion}). On the Lie algebra level,
this yields the semi-direct sum $\ah=\ah^{red}\psumma\ah^{nil}$.
The following assertion is a transfer to the double setting of
some properties of $\tri$-triples. Its proof
is a combination of arguments used in \cite[\S 2]{vi79} and
\cite[sect.\,1]{vitya}.
\begin{s}{Theorem}  \label{charact}
{\sf 1}. Each nilpotent pair has a characteristic; \\
{\sf 2}. If $\hb$ and $\hb'$ are two characteristic, then there exist $u\in
Z_G(\eb)^{nil}$ such that $u{\cdot}\hb=\hb'$;\\
{\sf 3}. If $\hb$ is any characteristic, then the eigenvalues of $\ad h_1$ and
$\ad h_2$ are rational.
\end{s}\begin{proof*}
1. Consider the algebraic group
$N:=\{g\in G\mid g{\cdot}e_i\in \Bbbk e_i, i=1,2\}$.
In view of condition \ref{nil-p}(ii), we have the exact sequence
\[ \{1\}\rar Z_G(\eb)\rar N\stackrel{\tau}{\rar} (\Bbbk^*)^2\rar \{1\} \ .
\]
(If $g{\cdot}e_i=t_ie_i$, then $\tau(g)=(t_1,t_2)$.)
Since $(\Bbbk^*)^2$ is reductive and Abelian, both $N^{nil}$ and
$[N,N]$ lie in $\Ker\tau$. Hence $\tau$ induces a surjective homomorphism
$(N^{red}/[N^{red},N^{red}])^o\rar (\Bbbk^*)^2$.
By a standard property of diagonalizable groups, this means that there
exists a 2-dimensional torus $C\subset N^o$ such  that $\vert Z_G(\eb)\cap
C\vert<\infty$ and $\tau\vert_C$ is onto. On the Lie algebra level, this
yields $\n^{nil}=\z_\g(\eb)^{nil}$ and $\n^{red}=\z_\g(\eb)^{red}\oplus\ce$.
The restriction of the Killing form to either of the reductive subalgebras
$\n^{red}$ and $\z_\g(\eb)^{red}$ is non-degenerate \cite[ch.\,4~\S 1.1]{vion}.
Hence $C$ can be chosen
so that $\ce\perp \z_\g(\eb)^{red}$. As $\n^{nil} \perp\n$, we obtain
$\ce\perp\z_\g(\eb)$. Restricting the differential of $\tau$ to $\ce$ yields
an isomorphism $d\tau:\ce\isom \Bbbk^2$. Define $h_1,h_2\in\ce$ by
$d\tau(h_1)=(1,0)$, $d\tau(h_2)=(0,1)$. Then $h_1,h_2$ have rational
eigenvalues and satisfy \re{comrel}.
\par
2. Suppose $\hb'=(h'_1,h'_2)$ is another characteristic of $\eb$.
Since $h'_i-h_i\in \z_\g(\eb)$ and $h'_i-h_i \perp \z_\g(\eb)$, we obtain
$h'_i-h_i\in \z_\g(\eb)^{nil}$ $(i=1,2)$. Thus $\ce$ and $\ce':=\langle
h'_1,h'_2\rangle$ are two maximal diagonalizable subalgebras of the solvable
algebra
$\ce\psumma\z_\g(\eb)^{nil}$. By the standard conjugacy theorem
(see \cite[ch.3~\S 2]{vion}), there exists
$u\in Z_G(\eb)^{nil}$ such that $u{\cdot}\ce=\ce'$. It then follows from
\re{comrel} that $u{\cdot}h_i=h'_i$ $(i=1,2)$.
\par
3. Because the characteristic constructed in the first part of the proof had
rational eigenvalues, we conclude by the second part.
\end{proof*}%
\begin{s}{Corollary} Suppose $\eb$ is a nilpotent pair and $\hb$ is a
semisimple pair satisfying \re{comrel}. If $\z_\g(\hb)\cap\z_\g(\eb)=0$,
then $\hb$ is a characteristic of $\eb$ and $\z_\g(\eb)$ contains no
semisimple elements.
\end{s}\begin{proof}
Let $\g=\z_\g(\hb)\oplus\me$ be the $Z_G(\hb)$-stable decomposition. Then
$\z_\g(\eb)\subset\me$ and therefore $\z_\g(\eb)$ is orthogonal to
$h_1,h_2$. That is, $\hb$ is a characteristic of $\eb$. From the proof of
\re{charact}, it follows that $\z_\g(\hb)$ contains a reductive Levi
subalgebra of $\z_\g(\eb)$. Thus, $\z_\g(\eb)^{red}=0$.
\end{proof}%
This theorem shows that nilpotent pairs are ``sufficiently good" double
analogues of nilpotent elements.
Willing to extend the classical theory of nilpotent orbits and
$\tri$-triples to the double setting, it is worth to
recall Dynkin's approach in 1952. Any 3-dimensional simple algebra $\ah$ has
the unique basis $\{e,\tih,f\}$ such that  $[\tih,e]=2e,[\tih,f]=-2f,[e,f]=\tih$.
The semisimple element $\tih$ is called the defining vector or characteristic
of $\ah$. Considering
$G$-orbits (conjugacy classes) of 3-dimensional simple subalgebras in $\g$,
Dynkin first proved that two such subalgebras are $G$-conjugate if and only if
their characteristics are \cite[Th.\,8.1]{EBD}.
Second, he proved that if
$\tih_+$ is the dominant representative in $G{\cdot}\tih$
(relative to a fixed set of simple roots $\ap_1,\dots,\ap_n$), then
$\ap_i(\tih_+)\in\{0,1,2\}$ \cite[Th.\,8.3]{EBD}.
This readily yields finiteness for the number of
$G$-orbits of characteristics and hence of nilpotent orbits.
According to the modern terminology, $\{e,\tih,f\}$ is called an $\tri$-triple
and $\tih$ is also called a characteristic of $e$.\\
It turns out that the second part of the above program has a counterpart in the
double setting. Let us fix some relevant notation and terminology.
Suppose $\hb$ is a characteristic of $\eb$. We shall consider the
bi-grading of $\g$ determined by $\hb$: \quad
$\g=\bigoplus_{p,q}\g_{p,q}$, where $\g_{p,q}=\{x\in\g\mid [h_1,x]=px,
[h_2,x]=qx\}$ and
$(p,q)$ ranges over a finite subset 
of ${\Bbb Q}\times{\Bbb Q}$ containing $(0,0)$, $(1,0)$, $(0,1)$. The pairs
$(p,q)$ with $\g_{p,q}\ne 0$ 
will be referred to as the {\it eigenvalues\/} of $\hb$ in $\g$. An
eigenvalue $(p,q)$ is said to be {\it integral}, if $p,q\in\Bbb Z$. Otherwise
it is called {\it fractional}. The same terminology is used for
the corresponding eigenspaces.
\begin{s}{Theorem}   \label{finita}
There exist finitely many $G$-orbits of characteristics of nilpotent pairs.
\end{s}\begin{proof}
1. Arguing by  induction on $\dim\g+\rk\g$, we first note that the claim
is true for $\tri$. \par
2. Assuming that $\hb$ has fractional eigenvalues, one may replace $\g$ by the
smaller semisimple subalgebra $[{\tilde\g},{\tilde\g}]$, where
${\tilde\g}=\oplus_{p,q\in {\bbb Z}}\g_{p,q}$. Indeed, $\tilde\g$ is reductive
and $e_1,e_2,h_1,h_2\in
{\tilde\g}$. Then, being nilpotent, $e_1,e_2$ belong to
$[{\tilde\g},{\tilde\g}]$. Since $h_1,h_2$ are orthogonal to $\z_\g(\eb)\cap
{\tilde\g}$ and the latter contains the centre of $\tilde\g$ (if any),
we have $h_1,h_2\in [{\tilde\g},{\tilde\g}]$. Thus, $\hb$ is a
characteristic of $\eb$ relative to $[{\tilde\g},{\tilde\g}]$ and, by the
inductive assumption, $[{\tilde\g},{\tilde\g}]$ contains finitely many
$[\tilde G,\tilde G]$-orbits of characteristics. Clearly these orbits generate
finitely many $G$-orbits in $\g\times\g$.
\par
3. Assume now that $\g={\tilde\g}$, i.e., the eigenvalues of
$\hb$ are integral. Fix a Cartan subalgebra $\te$ of $\g$ such that
$\{h_1,h_2\}\subset\te$. Then $\te\subset\g_{0,0}$.
Choose a set of simple roots $\Pi$
with respect to $\te$ so that $h_2+\esi h_1$ is dominant for all
sufficiently small positive $\esi\in \Bbb Q$.
For all $\ap\in\Pi$, we then have $\ap(h_2)\ge 0$ and if $\ap(h_2)=0$,
then $\ap(h_1)\ge 0$.
\par
4. Assume that $\beta(h_2)\ge 2$ for some $\beta\in \Pi$. Then we may
just throw it away! That is, consider $\Pi'=\Pi\setminus\{\beta\}$ and
the corresponding Levi subalgebra $\g'\subset\g$. The
constraint on $\beta$ implies that  $\g'_{p,q}=\g_{p,q}$ for $(p,q)\in
\{(0,0),(1,0),(0,1)\}$. Hence $e_1,e_2,h_1,h_2\in\g'$ and, as in part 2,
we even have $e_1,e_2,h_1,h_2\in [\g',\g']$. Thus, we may apply the inductive
assumption to $[\g',\g']$.
\par
5. Certainly, the constraint that $\ap(h_2)\in\{0,1\}$ for all $\ap\in\Pi$
leaves no much room and
one gets only finitely many possibilities for $h_2$. By symmetry, the same
argument applies to $h_1$. It follows that, up to conjugacy, there are
finitely many possibilities for $\hb$. This completes the proof.
\end{proof}%
Instead of appealing to the symmetry, one can exploit another argument in the
last part of the proof, keeping the same choice of $\Pi$. This has an
advantage of giving a more precise information about numerical labels
$\ap(h_i)$ ($i=1,2$). To this end, recall that a subalgebra
$\es\subset\g$ is called {\it regular},
if its normalizer $\n_\g(\es)$ contains a Cartan subalgebra.
Set $\el_i=\z_\g(h_i)$, $i=1,2$.
\begin{s}{Theorem}  \label{labels}
Let  $\eb$ be a nilpotent pair with a characteristic $\hb$. Suppose
$\{e_1,e_2\}$ is not contained in a proper regular semisimple subalgebra
of $\g$. Fix a Cartan subalgebra $\te$ containing
$\{h_1,h_2\}$. Then there exists a set of simple roots $\Pi$ relative
to $\te$ such that \par
{\rm (i)} $\ap(h_2)\in\{0,1\}$ for all $\ap\in\Pi$; \par
{\rm (ii)} If $\ap(h_2)=0$, then $\ap(h_1)\in\{0,1\}$; \par
{\rm (iii)} If $\ap(h_2)=1$, then $\ap(h_1)\in\{d,\dots,-1,0\}$, where $d$
is a (negative) constant depending only on $\el_2$.
\end{s}\begin{proof}
1. Inductive steps used in the proof of Theorem \ref{finita} provide us with
regular subalgebras in $\g$. Hence, under our hypothesis on $\eb$ and with the
same choice of $\te$ and $\Pi$, we already see that the eigenvalues of
$\hb$ must be
integral, $\ap(h_2)\in\{0,1\}$, and $\ap(h_1)\ge 0$, if $\ap(h_2)=0$. \\
2. If either $\ap(h_1)\ge 2$ or $\ap(h_1)=\ap(h_2)=1$ for some $\ap\in\Pi$,
one can again, as in
the proof of \re{finita}, throw away this $\ap$ and get a regular semisimple
subalgebra containing $e_1,e_2$. Thus, this cannot occur. \\
3. It remains to obtain the lower bound for $\ap(h_1)$, if $\ap(h_2)=1$.
Notice that $h_1$ induces in $\el_2$ the  $\Bbb Z$-grading
$\el_2=\oplus_{p\in{\Bbb Z}}\g_{p,0}$ and that $\Pi_2:=\{\ap\in\Pi\mid
\ap(h_2)=0\}$ is a set of simple roots for $\el_2$.
We also know that $\ap(h_1)\le 1$ for $\ap\in\Pi_2$.
Look what is happening on the next level. Suppose $\mu\in\Pi$ is such that
$\mu(h_2)=1$, and let $e_\mu$ be a nonzero root vector.
Then $e_\mu$ is a {\sl lowest\/} weight vector of an irreducible $\el_2$-module.
Denote this module by $\V_\mu$. We have $\V_\mu\subset\oplus_{p\in {\Bbb Z}}
\g_{p,1}$. Let $\lb$ be the root of $\g$ corresponding to the highest
weight of $\V_\mu$. Then, more precisely,
\[
\V_\mu\subset\bigoplus_{p=\mu(h_1)}^{\lb(h_1)}\g_{p,1} \ .
\]
Let us say that $\mbox{T}(\V_\mu):=
\lb(h_1)-\mu(h_1)$ is the height of $\Bbb Z$-grading on
$\V_\mu$. If $\lb(h_1)<0$, then $e_2\not\in \V_\mu$ and, as above, we
could drop the simple root $\mu$. Under our assumption, this is however
impossible and we must have $\lb(h_1)\ge 0$. It is thus enough to give an
upper bound on $\lb(h_1)-\mu(h_1)$ for any such $\mu$. The problem can be
stated as follows: \\[1ex]
{\sl The algebra $\el_2$ and an $\el_2$-module $\V_\mu$ are compatibly
$\Bbb Z$-graded; the
grade of each simple root in $\el_2$ is either 0 or 1. Give an upper bound
on the height of the $\Bbb Z$-grading on $\V_\mu$.} \\[1ex]
The answer was essentially given by Dynkin, who considered the height of
representation in case, where all simple roots are of grade 1
\cite[Suppl.,~\S 2\,$n^o$\,12]{max}. (Actually, Dynkin considered not
gradings but the associated partitions of the weights of a representation into
`layers', the height being \{the number of layers\}$-1$.)
If some of the simple roots have grade 0,
then the height can only decrease. Hence the bound given by Dynkin applies
in our situation as well. Being adapted to our setting, it reads
$\mbox{T}(\V_\mu)\le -2(\mu\vert\nu)$, where $\nu$ is the sum of all positive
coroots of $\el_2$ and $(\cdot\vert\cdot)$ is a Weyl group invariant inner
product on
$\te^*$. Whence $\mu(h_1)\ge 2(\mu\vert\nu)$, and we may
take $\displaystyle d=\min_{\mu\in\Pi\setminus\Pi_2} 2(\mu\vert\nu)$.
\end{proof}%
{\bf Remark.} In section 4, we describe some classes of nilpotent pairs
that do not lie in proper regular semisimple subalgebras and give a better
expression for $d$ in those cases.
\\[.6ex]
Unfortunately, attempts to extend further the Morozov-Jacobson theorem to
arbitrary nilpotent pairs fail: It is impossible in general to introduce
an ``opposite" nilpotent pair $\fb$. And what is worse, the number of
$G$-orbits of all nilpotent pairs is infinite (see \cite[5.5]{vitya}).
Comparing with Theorem \ref{finita} implies that there should exist infinite
families of $G$-orbits of nilpotent pairs with the same characteristic.
Actually, it is not
hard to realize that nilpotent pairs described in [loc.\,cit] give examples
of such families.
Therefore one can address the following somewhat vague problem:
\\[0.8ex]
\hspace*{1pt}\refstepcounter{equation} {\bf (\theequation)} \quad
\label{problem}
\parbox{390pt}{%
{\it Find a natural subclass of all nilpotent
pairs that consists of finitely many orbits, admits a rich structure theory,
and includes all interesting examples. }}
\\[1.2ex]
Let us describe available `interesting' examples.
We first indicate the most restrictive case, where the opposite $\fb$ exists
and the theory becomes entirely parallel to the ordinary one.
By the invariance of the Killing form, $\z_\g(\eb)^{\perp}=[\g,e_1]+[\g,e_2]$.
That is, $h_i\in \Ima(\ad e_1)+\Ima(\ad e_2)$, if $\hb$ is a characteristic
of $\eb$.
\begin{s}{Proposition}  \label{rectan}
The following conditions are equivalent: \par
{\sf 1.} $h_1\in\Ima(\ad e_1)$; \par
{\sf 2.} $h_2\in\Ima(\ad e_2)$; \par
{\sf 3.} there exist commuting $\tri$-triples $\{e_1,2h_1,f_1\}$ and
$\{e_2,2h_2,f_2\}$.
\end{s}\begin{proof}
Since (3) implies (1) and (2), it suffices to demonstrate that
(1)$\Rightarrow$(3). \\
Suppose $h_1=[e_1,f]$. Then replacing $f$ by its projection to the
$(-1)$-eigenspace of $\ad h_1$, we obtain the $\tri$-triple
$\{e_1,2h_1,f_1\}$. Further, $0=[e_2,h_1]=[e_1,[e_2,f_1]]$. Hence
$[e_2,f_1]$ lies in the $(-1)$-eigenspace of $\ad h_1$ in $\z_\g(e_1)$.
As the latter is positively graded, this forces $[e_2,f_1]=0$.
Similarly, $[h_2,f_1]=0$. Thus $e_2,h_2$ belong to the reductive
subalgebra $\ka_1:=\z_\g(e_1,h_1,f_1)$. Since the restriction of the Killing
form to $\ka_1$ is non-degenerate, it follows from the condition
$h_2\perp\z_\g(\eb)$ that $h_2$ is orthogonal to $\z_{\ka_1}(e_2)$ and hence
$h_2\in [\ka_1,e_2]$. Thus, $e_2$ and $2h_2$ can be included into the
$\tri$-triple inside of $\ka_1$.
\end{proof}%
{\bf Remark.} Each member of
a nilpotent pair has a characteristic in its own right. To avoid
a confusion between characteristics of nilpotent pairs and
of nilpotent elements, the latter are always marked with `tilde'.
We shall say that $\eb$ is {\it rectangular} whenever the equivalent
conditions of \re{rectan} hold. Roughly speaking, the proposition
claims that $\eb$ is rectangular if and only if $2h_i=\tih_i$ ($i=1,2$).
\\[1ex]
The number of $G$-orbits of rectangular pairs is evidently finite, and
various assertions for such pairs can immediately be derived on the base of the
ordinary theory (see \cite{wir},\,\cite[sect.\,3]{pCM}). Yet, there exist a
plenty of interesting not necessarily rectangular pairs: {\it principal\/}
\cite{vitya} and {\it almost principal\/} \cite{pCM} nilpotent pairs.
A nilpotent pair $\eb$ is called principal
(resp. almost principal) if
$\dim\z_\g(\eb)=\rk\g$ (resp. $\rk\g+1$). Dealing with these pairs, we shall
usually omit the adjective `nilpotent'.
Finiteness of the number of $G$-orbits of principal pairs was proved in
\cite[3.9]{vitya}. Below, we show this holds for the almost principal pairs.
A variety of other results obtained for principal and almost principal pairs
allows us to treat them as `very interesting'.
It is therefore natural to require that a subclass we are searching for
would include all rectangular, principal, and almost principal pairs.

\sekt{Wonderful nilpotent pairs\nopagebreak}%
In this section we present our solution to problem \ref{problem}.
The idea to consider filtrations and limits associated with nilpotent
elements is due to R.\,Brylinski \cite{ranee}. Then V.\,Ginzburg realized
that this nicely works in case of nilpotent pairs.
\\[1ex]
Let $\eb$ be a nilpotent pair in $\g$ and $\hb$ a characteristic of it.
Put ${\Bbb P}:=\{0,1,2,\dots\}$.
Making use of $\eb$, one may define three filtrations
for any subspace $M\subset\g$: \par
$\bullet$ $e_1$-filtration:
$M(i,\ast)=\{x\in M\mid (\ad e_1)^{i+1}x=0\}$, $i\in\Bbb P$\ ;
\par
$\bullet$ $e_2$-filtration:
$M(\ast,j)=\{x\in M\mid (\ad e_2)^{j+1}x=0\}$, $j\in\Bbb P$\ ;
\par
$\bullet$  $\eb$-filtration: $M(i,j)=M(i,\ast)\cap M(\ast,j)$. \\
The corresponding limits are defined by the formulas: \par
lim$_{e_1} M=\sum_{i\in {\Bbb P}} (\ad e_1)^i M(i,\ast)\subset\g$, \par
lim$_{e_2} M=\sum_{j\in {\Bbb P}} (\ad e_2)^j M(\ast,j)\subset\g$, \par
lim$_\eb M=
\sum_{i,j\in {\Bbb P}} (\ad e_1)^i(\ad e_2)^j M(i,j)\subset\g$.
\\[1ex]
Notice that the order of $\ad e_1$ and $\ad e_2$ in the last line is
immaterial, since $[e_1,e_2]=0$.
It follows that $\lim_{e_i}M\subset\z_\g(e_i)$ and $\lim_\eb M\subset
\z_\g(\eb)$.
It is easily seen that in all three cases $\dim\,(\mbox{lim\,}_\bullet M)
\le\dim M$ and the equality is equivalent to the fact that the sum
in the definition of corresponding limit is actually direct.
(To check this for the $\eb$-limit, one should use the equality \quad
$\dim(\ad e_1)^i(\ad e_2)^j M(i,j)=$\\
$\dim M(i,j){-}\dim M(i{-}1,j){-}
\dim M(i,j{-}1){+}\dim M(i{-}1,j{-}1)$.) \\
We shall repeatedly use the following sufficient condition for this to happen.
In case of ordinary filtrations this was observed in \cite{ranee}.
Recall that $\el_i=\z_\g(h_i)$.
\begin{s}{Lemma}   \label{ravno} \\
{\sf 1.} If $M\subset \el_i$, then $\dim (\lim_{e_i}M)=
\dim M$ ($i=1,2$);\\
{\sf 2.} If $M\subset \z_\g(\hb)$,
then $\dim\,(\lim_{\eb}M)=\dim M$.
\end{s}\begin{proof}
Relations \ref{comrel} show that different summands in definition of
all limits belong to different eigenspaces (relative to $h_i$ and
$\hb$ respectively).
\end{proof}%
If $A=\oplus_i A_i$ is a $\Bbb Q$-graded object and $M\subset\Bbb Q$, then
$A_M:=\oplus_{i\in M}A_i$. We shall use this for $M\in\{{\Bbb N},
{\Bbb P},{\Bbb Z}\}$.
Similar notation is used in the bi-graded situation. This will be applied to
various subspaces of $\g$, the bi-grading being determined by $\hb$.
For instance, $\z_\g(\eb)$ is ${\Bbb Q}{\times}\Bbb Q$-graded and we may
consider $\z_\g(\eb)_{{\Bbb P},{\Bbb P}}$. On the other hand,
$\z_\g(e_2,h_1)$ lies in $\el_1$ and hence is only $\Bbb Q$-graded.
Therefore $\z_\g(e_2,h_1)_{\Bbb P}$ is defined. It is meant here that
$\z_\g(e_2,h_1)_j$ is the subspace sitting in $\g_{0,j}=(\el_1)_j$.
\par
Set $\h:=\z_\g(\hb)=\g_{0,0}=\el_1\cap\el_2$ and consider the above
filtrations and limits for $M=\h$.
It follows immediately form the definitions of the limits and Eq.\,\re{comrel}
that \\
\refstepcounter{equation} {\bf (\theequation)} \quad
\label{inclu}
\parbox{390pt}{%
\hspace{2.5cm}$\begin{array}{l}
\lim_{e_1}\h\subset\z_\g(e_1,h_2)_{\Bbb P}=\z_{\el_2}(e_1)_{\Bbb P} \\
\lim_{e_2}\h\subset\z_\g(e_2,h_1)_{\Bbb P}=\z_{\el_1}(e_2)_{\Bbb P} \\
\lim_{\eb}\h\subset\z_\g(\eb)_{{\Bbb P},{\Bbb P}} \ .
\end{array}$
}
\\[1.2ex]
The next equalities follow from the fact that, by Lemma \ref{ravno}, all
three spaces have the same dimension and the first two lie in the last one:
\beq  \label{sovpad}
\textstyle \lim_{e_1}(\lim_{e_2}\h)=\lim_{e_2}(\lim_{e_1}\h)=\lim_{\eb}\h \ .
\eeq
Now we are ready to introduce our main object.
\begin{rem}{Definition}  \label{wonder}
A nilpotent pair $\eb$ is called {\it wonderful\/} in $\g$, if
$\lim_\eb\h=\z_\g(\eb)_{{\Bbb P},{\Bbb P}}$;
or, in a more explicit form, if
\begin{center}
$(\ad e_1)^i(\ad e_2)^j\h(i,j)=\z_\g(\eb)_{i,j}$ for all $i,j\in\Bbb P$.
\end{center}
\end{rem}%
Our primary goal is to describe properties of such pairs.
For the sake of completeness, we start with an easy result.
\begin{s}{Lemma}  \label{wond0}
Let $\es$ be a semisimple  subalgebra of $\g$ and $\eb\in\es\times\es$ a
pair of commuting nilpotent elements. Then \par
{\rm (i)} $\eb$ is nilpotent in $\es$ if and only if it is nilpotent in $\g$;
\par
{\rm (ii)} if $\eb$ is wonderful in $\g$, then it is wonderful in $\es$.
\end{s}\begin{proof*}
(i) The `only if' part is obvious. \\
Suppose $\eb$ is nilpotent in $\g$ and $(h_1,h_2)\in\g\times\g$ satisfies
Eq.\,\re{comrel}. Let $\g=\es\oplus{\frak m}$ be an $\es$-invariant
decomposition and $h_i=h_i^{(\es)}+h_i^{(\frak m)}$, $i=1,2$. Then
$(h_1^{(\es)},h_2^{(\es)})$ also satisfy Eq.\,\re{comrel}. This shows that
$N_S:=\{g\in S\mid g{\cdot}e_i\in\Bbbk e_i\}$ contains sufficiently many
elements to ensure that \ref{nil-p}(ii) holds for $S$ in place of $G$. \\
(ii) Let $\eb$ be wonderful in $\g$. It follows from (i) and Theorem
\ref{charact} that there exists a characteristic $\hb$ lying in
$\es\times\es$. The following is obvious.
\end{proof*}%
\begin{s}{Proposition}   \label{wond1}
Suppose $\eb$ is wonderful. Then \\
{\sf 1.} $\lim_{e_1}\h=\z_\g(e_1,h_2)_{\Bbb P}$ and
$\lim_{e_2}\h=\z_\g(e_2,h_1)_{\Bbb P}$; \\
{\sf 2.} $\lim_{e_1}\z_\g(e_2,h_1,h_2)=\z_\g(e_2,e_1,h_2)_{\Bbb P}$
and $\lim_{e_2}\z_\g(e_1,h_1,h_2)=\z_\g(e_1,e_2,h_1)_{\Bbb P}$.
\end{s}\begin{proof}
By symmetry, it is enough to prove one equality in each item.\\
1. Using \re{inclu} and \re{sovpad} we obtain \par
$\z_\g(\eb)_{{\Bbb P},{\Bbb P}}=\lim_\eb\h=
\lim_{e_2}(\lim_{e_1}\h)\subset\lim_{e_2}(\z_\g(e_1,h_2)_{\Bbb P})
\subset\z_\g(\eb)_{{\Bbb P},{\Bbb P}}$.\\
2. Applying the formula in Definition \ref{wonder} with $j=0$ gives
$(\ad e_1)^i\h(i,0)=\z_\g(\eb)_{i,0}$ for all $i\in\Bbb P$. Then summation
over $i$ yields the first formula.
\end{proof}%
The following simple lemma about $\Bbb Z$-graded Lie algebras appears to
be extremely useful in our situation.
\begin{s}{Lemma}   \label{useful}
Let $\ah=\oplus_{i\in {\Bbb Z}}\ah_i$ be a $\Bbb Z$-graded Lie algebra.
Suppose there exists $e\in\ah_1$ such that $\dim\z_\ah(e)_{\Bbb P}=\dim\ah_0$.
Then \par
{\sf 1.} $\lim_e\ah_0=\z_\ah(e)_{\Bbb P}$ and $[\ah_i,e]=\ah_{i+1}$ for all
$i\in\Bbb P$; \par
{\sf 2.} If $\ah$ is reductive, then $\z_\ah(e)=\z_\ah(e)_{\Bbb P}$ and
$e$ is Richardson in the nilpotent radical of the parabolic subalgebra
$\ah_{\ge 0}:=\oplus_{i\ge 0}\ah_i$.
\end{s}\begin{proof}
1. The first equality follows from the relations $\lim_e\ah_0\subset
\z_\ah(e)_{\Bbb P}$ and $\dim\ah_0=\dim(\lim_e\ah_0)$. The hypothesis on
$e$  also implies that the kernel of the map
$(\ad e)_{\ge 0}: \ah_{\ge 0}\rar \ah_{\ge 1}$ is of dimension
$\dim\ah_0$. Thus $(\ad e)_{\ge 0}$ is onto. \par
2. Using an invariant nondegenerate bilinear form on $\ah$ and
surjectivity of $(\ad e)_{\ge 0}$, one obtains $(\ad e)_{<0}:
\ah_{\le -1}\rar \ah_{\le 0}$ is injective. Hence $\z_\ah(e)$ is concentrated
in nonnegative degrees. The second claim is just a reformulation of the fact
that $(\ad e)_{\ge 0}$ is onto.
\end{proof}%
{\bf Remark.} If $\ah$ is reductive, the situation looks very much as if
$e$ were an even nilpotent element and the $\Bbb Z$-grading in question arose
from an $\tri$-triple containing $e$. This is not however always the case.
For instance, consider $\ah={\frak sl}_n$.  Here the $\Bbb Z$-grading
with characteristic $(10\dots0)$ satisfies the hypothesis in \re{useful},
but it does not correspond to an $\tri$-triple. Moreover, the corresponding
nilpotent element is not even.
\begin{s}{Proposition}  \label{wond2} If $\eb$ is wonderful, then \\
{\sf 1.} $\lim_{e_1}\h=\z_\g(e_1,h_2)_{\Bbb Z}=\z_\g(e_1,h_2)_{\Bbb P}$ and
$[\g_{i,0},e_1]=\g_{i+1,0}$ for all $i\in\Bbb P$, \\
\hspace*{10pt}
$\lim_{e_2}\h=\z_\g(e_2,h_1)_{\Bbb Z}=\z_\g(e_2,h_1)_{\Bbb P}$
and $[\g_{0,j},e_2]=\g_{0,j+1}$ for all $j\in\Bbb P$;\\
{\sf 2.} $[\z_\g(e_2,h_2)_i,e_1]=\z_\g(e_2,h_2)_{i+1}$
for all $i\in\Bbb P$, \\
\hspace*{10pt} $[\z_\g(e_1,h_1)_j,e_2]=\z_\g(e_1,h_1)_{j+1}$ for all $j\in\Bbb P$.
\end{s}\begin{proof*}
1. In view of Proposition \ref{wond1}(1), the previous Lemma applies
to reductive Lie algebras $(\el_1)_{\Bbb Z}$ and $(\el_2)_{\Bbb Z}$. \\
2. In view of Proposition \ref{wond1}(2), the previous Lemma applies
to Lie algebras $\z_\g(e_1,h_1)_{\Bbb Z}$ and $\z_\g(e_2,h_2)_{\Bbb Z}$.
\end{proof*}%
\begin{s}{Theorem}    \label{wond3}
Let $\eb$ and $\eb'$ be two wonderful pairs with the same characteristic $\hb$.
Then there exists $s\in H:=Z_G(\hb)$ such that $s{\cdot}e_i=e'_i$, $i=1,2$.
\end{s}\begin{proof}
We have $e_1,e'_1\in\g_{1,0}\subset (\el_2)_{\Bbb Z}$ and
$e_2,e'_2\in\g_{0,1}\subset (\el_1)_{\Bbb Z}$. By proposition \ref{wond2}(1)
with $i=0$, $[\g_{0,0},e_1]=\g_{1,0}$. This means that $e_1$
(and also $e'_1$) lies in the dense $H$-orbit in $\g_{1,0}$. Hence we may
assume that $e_1=e'_1$. Let $H_{e_1}$ denote the stabilizer of $e_1$ in $H$.
Then $\z_\g(e_1,h_1,h_2)=\z_\g(e_1,h_2)_0$ is Lie algebra of $H_{e_1}$.
By proposition \ref{wond2}(2)
with $j=0$, $[\z_\g(e_1,h_1)_0,e_2]=\z_\g(e_1,h_1)_1$. This means that $e_2$
(and hence $e'_2$) lies in the dense $H_{e_1}$\!-orbit in $\z_\g(e_1,h_1)_1$.
Thus we can make $e_2=e'_2$.
\end{proof}%
Next claim is a straightforward corollary of Theorems \ref{finita} and
\ref{wond3}:
\\[1ex]
\refstepcounter{equation} {\bf (\theequation)} \quad
\label{wondfin}
{\sl There exist finitely many $G$-orbits of wonderful pairs.}
\\[1.2ex]
Notice that the definition of wonderful pairs concerns only properties
of integral eigenspaces of $\hb$.
There is therefore no harm in assuming
that the eigenvalues of $\hb$ in $\g$ are integral. In this case, we say
that $\eb$ is {\it integral}.
\begin{s}{Lemma} \label{int} $\eb$ integral $\Leftrightarrow$
$\z_\g(\eb)=\z_\g(\eb)_{{\Bbb Z},{\Bbb Z}}$.
\end{s}\begin{proof}
If $(\ap,\beta)$ is a fractional eigenvalue of $\hb$ in $\g$, then, applying
nilpotent endomorphisms $\ad e_1$ and $\ad e_2$ to $\g_{\ap,\beta}$, we
eventually obtain an eigenvalue in $\z_\g(\eb)$ of the form $(\ap+m,\beta+n)$
with $m,n\in\Bbb P$.
\end{proof}%
In general, the passage from $\g$ to
$\g_{{\Bbb Z},{\Bbb Z}}$ can be considered as double analogue of taking
the {\it even} part of
the $\Bbb Z$-grading associated to an $\tri$-triple. (A discrepancy is
explained by the fact that, unlike Eq.\,\re{comrel}, the standard
normalization in $\tri$-triples is: $[h,e]=2e$.)
For future references, we record the following fact
which is a straightforward consequence of Lemma \ref{useful}(2) and
Proposition \ref{wond2}(1):
\begin{s}{Lemma}  \label{richardson}
Let $\eb$ be wonderful and integral. Then $e_1$ is Richardson in
$\el_2$ and $e_2$ is Richardson in $\el_1$.
\end{s}%
The definition of a wonderful pair says something about $\z_\g(\eb)$ in
the positive quadrant (of ${\Bbb Z}^2$-grading). This allows us to draw a
conclusion about the negative quadrant.
\begin{s}{Proposition}  \label{pusto3}
Let $\eb$ be wonderful and integral. Then $\z_\g(\eb)_{p,q}=0$ whenever
$p<0,q<0$.
\end{s}\begin{proof}
The argument used in the proof of Theorem 2.5(1) in \cite{pCM} applies here
verbatim. For convenience of the reader, we reproduce it. \\[.6ex]
Assume $\z_\g(\eb)_{p,q}$ is nonzero for $p_0=-p>0$ and
$q_0=-q>0$. It follows from the invariance of the Killing form on $\g$ that
$\z_\g(\eb)_{p,q}\ne 0$ if and only if
$\g_{p_0,q_0}\not\subset \mbox{Im\,}(\ad e_1)+\mbox{Im\,}(\ad e_2)$.
By definition, put ${\cal D}=
\g_{p_0,q_0}\setminus (\mbox{Im\,}(\ad e_1)+\mbox{Im\,}(\ad e_2))$. For each
$y\in\cal D$, consider the finite set $I_y=\{(k,l)\in ({\Bbb Z}_{\ge 0})^2\mid
(\ad e_1)^k(\ad e_2)^ly\ne 0\}$, with the lexicographic ordering.
This means $(k,l)\prec(k',l')$ $\Leftrightarrow$
$k<k'$ or $k=k'$ and $l<l'$. Denote by $m(I_y)$ the unique maximal element
in $I_y$. Let $y^*\in\cal D$ be an element such that
$(k_0,l_0):=m(I_{y^*})\preceq m(I_z)$ for all $z\in\cal D$. Then
$(\ad e_1)^{k_0}(\ad e_2)^{l_0}y^*$ is a nonzero element in $\z_\g(\eb)\cap
\g_{p_0+k_0,q_0+l_0}$. By Definition \ref{wonder},
there is $t\in\h(p_0+k_0,q_0+l_0)$
such that $(\ad e_1)^{p_0+k_0}(\ad e_2)^{q_0+l_0}t=
(\ad e_1)^{k_0}(\ad e_2)^{l_0}y^*$. Then $(\ad e_1)^{k_0}(\ad e_2)^{l_0}(y^*-
(\ad e_1)^{p_0}(\ad e_2)^{q_0}t)=0$. Since $p_0>0,q_0>0$, we have
$z^*=y^*-(\ad e_1)^{p_0}(\ad e_2)^{q_0}t$ is nonzero and belongs to $\cal D$.
However, $I_{z^*}\subset I_{y^*}\setminus\{(k_0,l_0)\}$. Therefore
$m(I_{z^*})< m(I_{y^*})$, which contradicts the choice of $y^*$.
Thus, the case $p<0,q<0 $ is impossible.
\end{proof}%
It appears that combining this proposition with information about
$\z_{\el_1}(e_2)$ and $\z_{\el_2}(e_1)$ contained in \ref{wond2}, one obtains
a characterization of the wonderful pairs.
\begin{s}{Theorem}   \label{xarak}
Let $\eb$ be integral. Then
\begin{center}
$\lim_\eb\h=\z_\g(\eb)_{{\Bbb P},{\Bbb P}}\quad \Longleftrightarrow
\quad \left\{
\begin{array}{lll}
\z_\g(\eb)_{p,q}=0 & \mbox{ for } & p<0,q<0, \\
\z_\g(e_1)_{p,0}=0 & \mbox{ for } & p<0, \\
\z_\g(e_2)_{0,q}=0 & \mbox{ for } & q<0 \quad .
\end{array}
\right.
$
\end{center}
\end{s}\begin{proof}
``$\Rightarrow$" \ This is Propositions \ref{pusto3} and \ref{wond2}(1). \\
``$\Leftarrow$"  \ Making use of the Killing form on $\g$, one can translate
these three conditions in ones about the positive quadrant. Namely, \par
1st:\quad $\g_{p,q}\subset \Ima(\ad e_1)+\Ima(\ad e_2)$ for $p>0,q>0$; \par
2nd:\quad $\g_{p,0}\subset \Ima(\ad e_1)$ for $p>0$; \par
3rd:\quad $\g_{0,q}\subset \Ima(\ad e_2)$ for $q>0$.
\\[1ex] These three together show that $\g_{0,0}(=\h)$, $e_1$, and
$e_2$ generate
$\g_{{\Bbb P},{\Bbb P}}$. In particular, applying $\ad e_1$ and $\ad e_2$ to
$\h$, we obtain the whole space $\z_\g(\eb)_{{\Bbb P},{\Bbb P}}$.
\end{proof}%
Of course, all statements about integral wonderful pairs can be reformulated
as ones about integral eigenspaces of arbitrary wonderful pairs.

\sekt{Classes of wonderful pairs\nopagebreak}%
First of all, notice that a rectangular nilpotent pair is wonderful.
This is an exercise for the reader. Next, the principal and almost
principal pairs are wonderful. This follows from \cite[sect.\,1]{vitya}
for the former and from \cite[2.3]{pCM} for the latter.
As we have now the general concept of a characteristic of a nilpotent pair,
it is worth to put these notions in a more general context.
\\[1ex]
Recall that the even nilpotent elements in $\g$ are characterized by the
following property: \par
Let $\{e,\tih,f\}$ be an $\tri$-triple. Then $e$ is even if and only if
$\dim\z_\g(e)=\dim\z_\g(\tih)$.
\\[1ex]
Translating this into the double setting, we shall say that a nilpotent
pair $\eb$ is \\
\hspace*{1cm} $\bullet$ {\it even\/} whenever $\dim\z_\g(\eb)=\dim\z_\g(\hb)$;
\\
\hspace*{1cm} $\bullet$ {\it almost even\/} whenever
$\dim\z_\g(\eb)=\dim\z_\g(\hb)+1$.
\\
Notice that the second assumption implies that $\eb$ is non-trivial, because
$\dim\z_\g(e)-\dim\z_\g(h)$ is always even for the sole nilpotent element.
Since $\rk\g\le \dim\z_\g(\hb)\le\dim\z_\g(\eb)$ and $\dim\z_\g(\hb)-
\rk\g$ is even, any principal pair is even and almost principal pair is almost
even.
Our terminology is partly justified by the following observation.
\begin{s}{Proposition}  \label{pr-even}
Let $\eb$ be a rectangular nilpotent pair. Then $\eb$ is even if and only if
both $e_1$ and $e_2$ are even nilpotent elements.
\end{s}\begin{proof}
Let $\{e_1,\tih_1,f_1\}$ and $\{e_2,\tih_2,f_2\}$ be commuting $\tri$-triples
and let $\g=\oplus_{i\in{\Bbb Z}}\g(i)$ be the grading
determined by $\tih_1$. Then $e_1\in\g(2)$ and $e_2\in\g(0)$. Set $\ka_1=
\z_\g(e_1,\tih_1,f_1)$. Recall 
that $(\tih_1,\tih_2)=2\hb$ in the rectangular case.
It is easily seen that $\z_\g(e_1)$ and
$\g(0)\oplus\g(1)$ are isomorphic $\ka_1$-modules. (In the notation
of sect.\,2, $\lim_{e_1}(\g(0)\oplus\g(1))=\z_\g(e_1)$.)
As $e_2\in\ka_1$, we obtain
\[
\dim\z_\g(e_1,e_2)=\dim\z_{\g(0)}(e_2)+\dim\z_{\g(1)}(e_2) \ .
\]
Since $\dim\z_\g(h_1,h_2)=\dim\z_{\g(0)}(h_2)\le\dim\z_{\g(0)}(e_2)$, we
see that $\eb$ is even if and only if $\g(1)=0$ (i.e. $e_1$ is even) and
$e_2$ is even in $\g(0)$. Then either by symmetry or by a direct
argument one concludes that $e_2$ is actually even in $\g$.
\end{proof}%
Using the notation of the previous proof, it is easy to state the similar
condition for $\eb$ being almost even.
\begin{s}{Proposition}  \label{pr-almeven}
A rectangular nilpotent pair $(e_1,e_2)$ is almost even if and only if
$e_2$ is even in $\g(0)$ and $\g(1)$ is an irreducible $\langle
e_2,\tih_2,f_2\rangle$-module.
\end{s}\begin{proof}
Left to the reader.
\end{proof}%
It turns out that many statements about principal and almost principal
pairs proved in \cite{vitya} and \cite{pCM} remain true, with essentially
the same proofs, for the even and almost
even pairs. Not trying to be exhaustive in this %
rewriting, we demonstrate several results. Recall that $\h=\z_\g(\hb)$.
\begin{s}{Theorem}  \label{even}
{\sf 1.} The following conditions are equivalent:\par
{\rm (i)} $\lim_\eb\h=\z_\g(\eb)$, \par
{\rm (ii)} $\z_\g(\eb)_{{\Bbb P},{\Bbb P}}=\z_\g(\eb)$;\\
{\sf 2.} Any even nilpotent pair is wonderful and integral.
\end{s}\begin{proof}
1. $(i)\Rightarrow (ii)$ -- We always have
$\lim_\eb\h\subset \z_\g(\eb)_{{\Bbb P},{\Bbb P}}\subset
\z_\g(\eb)$. \\
$(ii)\Rightarrow (i)$. As the centralizer of $\eb$ entirely lies in the
positive quadrant, we have $(\ad e_1)_{p,q}:\g_{p,q}\rar\g_{p+1,q}$
is injective for all $q$ and $p<0$. (Otherwise, applying $\ad e_2$ to
a nonzero element in the kernel we would eventually arrived at an element in
$\z_\g(\eb)_{p,q'}$ with $q'\ge q$.) Similarly,
$(\ad e_2)_{p,q}:\g_{p,q}\rar\g_{p,q+1}$
is injective for all $p$ and $q<0$. We thus have more than enough to apply
Theorem \ref{xarak} and conclude that $\lim_\eb\h=
\z_\g(\eb)_{{\Bbb P},{\Bbb P}}$. \\ \indent
2. This follows from the first part and from Lemma \ref{int}.
\end{proof}
The following sufficient condition will be helpful in our study of almost
even nilpotent pairs.
\begin{s}{Proposition}   \label{helpful}
Let $\eb$ be a nilpotent pair in $\g$. Suppose $\dim \V^{\langle
e_1,e_2\rangle}=1$ for a self-dual $\g$-module $\V$. Then $\eb$
is a rectangular principal nilpotent pair in $\g$.
\end{s}\begin{proof*} We first prove that $\eb$ has prescribed properties
as nilpotent pair in ${\frak sl}(\V)$ and then go down to $\g$.\\
1. Take a characteristic $\hb$ of $\eb$ and consider the corresponding
bi-grading
$\V=\oplus_{p,q\in{\Bbb Q}}\V_{p,q}$. Set $\Gamma=\{(p,q)\mid
\V_{p,q}\ne 0\}\subset {\Bbb Q}{\times}\Bbb Q$. Each nonempty coset
$((p',q')+({\Bbb Z}\times{\Bbb Z}))\cap\Gamma$ determines a subspace in $\V$
containing a $\langle e_1,e_2\rangle$-fixed vector (cf. \ref{int}).
Hence $\Gamma$ lies in a unique such coset.
For the same reason, $\dim\V_{p,q}=1$ for all $(p,q)\in\Gamma$, and
$\Gamma$ has a unique `northeast' corner, i.e.
\#$\{(p,q)\in\Gamma\mid (p{+}1,q)\not\in\Gamma\ \&
\ (p,q{+}1)\not\in\Gamma\}=1$. Let $(p_0,q_0)$ be this corner.
Since $\V$ is self-dual,
$\Gamma$ is centrally-symmetric. Whence $(-p_0,-q_0)$ is the unique
`southwest' corner of it. (Note that, although $p_0,q_0$ are not necessarily
integral, $(p_0,q_0)\in (-p_0,-q_0)+
({\Bbb Z}{\times}{\Bbb Z})$ implies $p_0,q_0\in\frac{1}{2}{\Bbb Z}$.)
It follows that $\Gamma$ lies inside of the
rectangle having opposite vertices $(p_0,q_0)$ and $(-p_0,-q_0)$.
It is however easy to see that the conditions $[e_1,e_2]=0$ and
$\dim \V^{\langle e_1,e_2\rangle}=1$ force that $\Gamma$ is ``equal" to
this rectangle, i.e., $\Gamma=\{(m,n)\mid
p_0-m\in {\Bbb Z}, q_0-n\in {\Bbb Z}, |m|\le p_0, |n|\le q_0\}$.
 \par
2. Since the eigenvalues of $\hb$ in $\V$ form a rectangle, $\eb$ is
a rectangular principal nilpotent pair in ${\frak sl}(\V)$. Indeed,
it is not hard to write a formula for the nilpotent operators $f_i$ such that
$[e_j,f_i]=\delta_{i,j}2h_j$ and $[h_j,f_i]=-\delta_{i,j}f_i$
($i,j\in\{1,2\}$).
Hence $\eb$ is rectangular. Set $\ah_i=\langle e_i,h_i,f_i\rangle$.
As ${\frak gl}(\V)\simeq\V\otimes\V^*$ and $\V$ is an irreducible
$\ah_1\times\ah_2$-module, an explicit computation shows
that $\dim {\frak gl}(\V)^{\langle e_1,e_2\rangle}=\dim\V=\rk{\frak gl}(\V)$
(cf. \cite{wir}).
Hence $\eb$ is principal in ${\frak sl}(\V)$. \par
3. Now, one has the following: $\g\subset {\frak sl}(\V)$ is semisimple,
$\{e_i,2h_i,f_i\}$ is an $\tri$-triple in ${\frak sl}(\V)$, and
$e_i,2h_i\in\g$. This easily implies that $f_i\in\g$ and $\eb$ is rectangular
in $\g$. \par
4. It follows from the general theory of principal nilpotent pairs
\cite[sect.\,1]{vitya} and is also easily
seen in our situation that $\z_{{\frak sl}(\V)}(\hb)$ is a Cartan subalgebra
(the set of diagonal traceless matrices).
Whence $\z_\g(\hb)$ is a Cartan subalgebra in
$\g$. Because the eigenvalues of $2h_i$ in $\g$ are even ($i=1,2$),
all irreducible $\ah_1{\times}\ah_2$-submodules in ${\frak sl}(\V)$,
and hence in $\g$, have zero weight. Thus $\dim\z_\g(\hb)=\dim\z_\g(\eb)$ and
$\eb$ is principal in $\g$.
\end{proof*}
\begin{rem}{Remarks}  \label{p0q0}
1. It follows from the proof that $\V$ is an irreducible
$\ah_1{\times}\ah_2$-module and, $(p_0,q_0)$ being the eigenvalue of
$\hb$ on $\V^{\langle e_1,e_2\rangle}$, $\dim\V=(2p_0+1)(2q_0+1)$. \\
2. Since $\V$ is assumed to be self-dual, $\g$ lies in either
${\frak sp}(\V)$ or ${\frak so}(\V)$. Obviously, $\g\hookrightarrow
{\frak so}(\V)$ if and only if $\V$ is an orthogonal
$\ah_1{\times}\ah_2$-module if and only if $p_0-q_0\in\Bbb Z$.
In any case, $\eb$ is principal in the respective classical Lie algebra.
Furthermore, it is not hard to list all
possibilities for such $(\g,\V)$, see \cite{wir2}.
\vskip-2ex
\end{rem}
\begin{s}{Theorem}  \label{almeven}
Suppose $\eb$ is almost even. Then $\z_\g(\eb)=\lim_\eb\h\oplus
\langle x\rangle$, where $x\in\g_{p,q}$. \\
{\sf 1}. There are two possibilities
for the eigenvalue of $x$:
either $p,q\in\Bbb Z$ and $pq<0$ or $p,q\in \frac{1}{2}{\Bbb Z}\setminus
{\Bbb Z}$ and $p>0,q>0$. In particular, $\eb$ is wonderful. \\
{\sf 2}. In case $(p,q)$ is fractional, we have \\ \indent
$\h$ is a Cartan subalgebra in $\g$, \\ \indent
$\eb$ is an almost principal rectangular pair,  \\ \indent
$\eb$ is principal in $\g_{{ \Bbb Z},{\Bbb Z}}$.
\end{s}\begin{proof}
Since $\lim_\eb\h$ is $\hb$-stable, the first equality follows for dimension
reason. \\
1. (a) Suppose $(p,q)$ is integral. That the case $p\ge 0,q\ge 0$ is impossible
follows from Theorem \ref{even}(1). This already means that $\eb$ is wonderful.
Consequently, results of sect.\,2 apply. By Theorem \ref{xarak}, neither
$p<0,q<0$ nor $pq=0$ can occur. We are thus left with the case $pq<0$.

(b) Suppose $(p,q)$ is fractional. Let $\g_{f\! r}$ be the sum of all
fractional eigenspaces. Then $\g_{f\! r}$ is an orthogonal
$\g_{{\Bbb Z},{\Bbb Z}}$-module, $\g=\g_{{\Bbb Z},{\Bbb Z}}\oplus\g_{f\! r}$,
and $\g_{f\! r}^{\langle e_1,e_2\rangle}=\langle x\rangle$. Applying
Proposition \ref{helpful}, we conclude that $p,q\in\frac{1}{2}{\Bbb Z}$,
$\eb$ is rectangular principal in $\g_{{\Bbb Z},{\Bbb Z}}$,
and $\z_\g(\hb)$ is a Cartan subalgebra of $\g$.
Moreover, it follows from the orthogonality that both $p,q$ must be
fractional (cf. Remark \ref{p0q0}).
Finally, $\dim\z_\g(\eb)=\rk\g_{{\Bbb Z},{\Bbb Z}}+1=\rk\g+1$, i.e.,
$\eb$ is almost principal in $\g$.
\end{proof}%
{\bf Remark.} Notice that in this case the decomposition
$\g=\g_{{\Bbb Z},{\Bbb Z}}\oplus\g_{f\! r}$ is a ${\Bbb Z}_2$-grading and the
corresponding involutory automorphism of $\g$ is inner.
\begin{rem}{Example}  \label{denom}
If $\eb$ is rectangular and $\hb$ is a characteristic  of it, then the
eigenvalues of $h_i$ belong to $\frac{1}{2}\Bbb Z$ ($i=1,2$). It was also
shown above that for the even nilpotent pairs the eigenvalues of $\hb$ are
integral.
There exist however wonderful pairs whose eigenvalues have arbitrarily large
denominators. An example can be constructed as follows. Let $\eb$ be a principal
nilpotent pair in ${\frak sl}_n$ corresponding to the partition $(n-1,1)$.
(See \cite[sect.\,5]{vitya} for the relationship between principal nilpotent
pairs in ${\frak sl}_n$ and partitions.)
Consider $\eb$ as nilpotent pair in ${\frak sp}_{2n}$, using the
natural inclusion ${\frak sl}_n\hookrightarrow {\frak sp}_{2n}$. Then $\eb$ is
a non-integral wonderful pair in ${\frak sp}_{2n}$ and eigenvalues of $\hb$
in ${\frak sp}_{2n}=\g$ have denominator `$n$'. For simplicity, take $n=3$.
Having chosen a Witt basis in $\Bbbk^6$ for the symplectic form,
we may take $e_1=v_{23}-v_{45}$ and
$e_2=v_{13}-v_{46}$, where $\{v_{ij}\}$ is the monomial basis in the space
of matrices. Then $h_1=\frac{1}{3}\mbox{diag}(-1,2,-1,1,-2,1)$,
$h_2=\frac{1}{3}\mbox{diag}(2,-1,-1,1,1,-2)$. Here $\dim\z_\g(\eb)=7$ and
the eigenvalues of $\hb$ on $\z_\g(\eb)$ are $(0,0),\,(1,0),\,(0,1),$
$(1/3,1/3),\,(2/3,2/3),\,(4/3,-2/3),\,(2/3,-4/3)$.
\end{rem}%
Thus, there is no universal bound for denominators of the eigenvalues of
$\hb$. Yet, one can give a (very rough) estimate for each $\g$, which actually
applies to characteristics of arbitrary nilpotent pairs. If $\es$ is
a semisimple Lie algebra, let $c(\es)$ denote the determinant of the Cartan
matrix of $\es$ or, equivalently, the order of the centre of the corresponding
simply-connected group.
\begin{s}{Lemma}  \label{lembo}
Let $h\in\es$ be a semisimple element. Suppose the eigenvalues of $\ad h$
are integral. For any finite-dimensional representation $\rho:\es\to
{\frak sl}(\V)$, then
the eigenvalues of $\rho(h)$ belong to $\frac{1}{c(\es)}\Bbb Z$.
\end{s}\begin{proof*}
The rows of the inverse of the Cartan matrix yield the expressions of the
fundamental weights through the simple roots.
\end{proof*}
\begin{s}{Proposition}  \label{estimate}
Let $\hb$ be a characteristic of a nilpotent pair in $\g$. Then the
denominators of the eigenvalues of $\ad h_1$, $\ad h_2$ do not exceed
$\displaystyle\max_{\es\subset\g}c(\es)$, where $\es$ ranges over all regular
semisimple subalgebras of $\g$.
\end{s}\begin{proof}
As in \re{almeven}, consider the decomposition $\g=\g_{{\Bbb Z},{\Bbb Z}}\oplus\g_{f\! r}$
and set $\es=[\g_{{\Bbb Z},{\Bbb Z}},\g_{{\Bbb Z},{\Bbb Z}}]$. Then $h_1,h_2\in\es$ and $\g_{f\! r}$ is
an $\es$-module. Now we conclude by the previous lemma.
\end{proof}%
However I think that for wonderful pairs a better estimate ought to exist.

\sekt{Characteristics for principal and almost principal pairs\nopagebreak}%
We give a version of Theorem \ref{labels} for principal and almost principal
integral pairs\footnote{in \cite{pCM}, we used the term ``pairs of
$\Bbb Z$-type" for integral pairs.}. When we shall give two references for a
property of such pairs, this means that the proof is found in \cite{vitya}
for principal pairs and in \cite{pCM} for almost principal ones.
\begin{s}{Lemma}  \label{notcontained}
Let $\eb$ be either a principal or an almost principal integral pair.
Then $\eb$
is not contained in a proper regular semisimple subalgebra of $\g$.
\end{s}\begin{proof}
(Cf. Remark after 4.4 in \cite{pCM}.)
Assume that $\eb$ is contained in a proper regular semisimple subalgebra
and let
$\tilde\g$ be a maximal one among them. It follows from
\cite[\S 5]{EBD} and V.\,Kac's description of periodic automorphisms of
$\g$ (see e.g. \cite[3.7]{t41})
that $\tilde\g$ is a fixed-point subalgebra of some {\it inner\/}
automorphism of $\g$ of finite order. Since $e_1,e_2\in\tilde\g$, this means
$Z_G(\eb)$ contains non-trivial semisimple elements. But $Z_G(\eb)$ is
connected unipotent for such $\eb$ (see \cite[3.6]{vitya} and
\cite[2.14]{pCM}).
\end{proof}%
By Ginzburg's result \cite[3.4]{vitya}, both $G{\cdot}e_1$ and $G{\cdot}e_2$
are Richardson orbits in $\g$, if $\eb$ is principal. But in the almost
principal case only one of them is Richardson \cite[2.9(ii),\,2.10]{pCM}.
In either case, if the orbit $G{\cdot}e_2$ is Richardson, a more
precise statement is:\\
Set $\displaystyle \p_2=\bigoplus_{i\in{\Bbb Z},j\in{\Bbb P}}\g_{i,j}=
\g_{\ast,{\Bbb P}}$. It is a parabolic subalgebra and $\el_2$ is a Levi
subalgebra in it. Then $e_2\in(\p_2)^{nil}=\g_{\ast,{\Bbb N}}$
and $[\p_2,e_2]=(\p_2)^{nil}$.
\\[.5ex]
For a reductive Lie algebra $\el$, let $\cox(\el)$ denote maximum among
the Coxeter numbers of the simple components of $\el$.
Recall that $\z_\g(\hb)$ is a Cartan subalgebra in the principal and almost
principal case.

\begin{s}{Theorem}  \label{pr-char}
Let $\eb$ be either a principal or an almost principal integral pair.
Choose the
set of simple roots $\Pi$ relative to $\z_\g(\hb)$ so that $h_2+\esi h_1$ is
strictly dominant for all sufficiently small $\esi\in\Bbb Q$. Then \par
{\rm (i)} $\ap(h_2)\in\{0,1\}$ for all $\ap\in\Pi$; \par
{\rm (ii)} If $\ap(h_2)=0$, then $\ap(h_1)=1$; \par
{\rm (iii)} If $\ap(h_2)=1$ and $G{\cdot}e_2$ is Richardson, then
$\ap(h_1)\in\{-\cox(\el_2){+}1,\dots,-1,0\}$.
\end{s}\begin{proof}
Using Theorem \ref{labels}, Lemma \ref{notcontained}, and the fact that
$\z_\g(\hb)$ is Cartan (hence the case $\ap(h_1)=\ap(h_2)=0$ is impossible),
one sees that we have to only prove that $\ap(h_1)\ge -\cox(\el_2){+}1$.
By \cite[sect.\,1]{vitya} and \cite[2.3]{pCM}, $e_1$ is regular nilpotent
in $\el_2$. It then follows from (ii) that the $\Bbb Z$-grading in $\el_2$
defined by $h_1$ is nothing but the standard grading associated with the
function $\ap\mapsto\hot(\ap)$ on the set of roots of $\el_2$ (i.e.,
$(\el_2)_i$ is the linear span of the root spaces such that the height
of the corresponding root of $\el_2$ is $i$). Therefore
$\min\{i\mid \g_{i,0}\ne 0\}=-\cox(\el_2){+}1$. If $e_2$ is Richardson in
$(\p_2)^{nil}$, then $[\g_{\ast,0},e_2]=\g_{\ast,1}$. Hence
$\min\{i\mid \g_{i,1}\ne 0\}\ge-\cox(\el_2){+}1$, which is exactly what we
need.
\end{proof}%
Recall that the exponents of a simple Lie algebra are the degrees of
fundamental polynomial invariants of the adjoint representation, reduced by 1.
The set of exponents of a reductive Lie algebra is the union of the exponents
of all simple components. Since the sum of the exponents is the dimension
of a maximal nilpotent subalgebra, the following is a generalization of
\cite[6.13]{vitya}.
\begin{s}{Theorem}  \label{exponents}
Let $\eb$ be either principal or almost principal and
let $(\ap_i,\beta_i)$ ($i=1,\dots,\rk\g$) be the eigenvalues of $\hb$ in
$\lim_\eb\h=\z_\g(\eb)_{{\Bbb P},{\Bbb P}}$.
Then $\{\ap_i\mid \beta_i\ne 0\}$ are the exponents of $\el_2$
and $\{\beta_i\mid \ap_i\ne 0\}$ are the exponents of $\el_1$.
\end{s}\begin{proof*}
1. First, assume that $\eb$ is either principal or almost principal integral.\\
Since $\eb$ is wonderful and integral in both cases,
the formulas in \re{wond2} become
simpler. In particular, $\lim_{e_1}\h=\z_\g(e_1,h_2)=\z_{\el_2}(e_1)$.
Using \re{sovpad}, we obtain $\lim_\eb\h=
\lim_{e_2}\z_{\el_2}(e_1)$. Because $[h_1,e_2]=0$,
this implies that the eigenvalues of $h_1$ in $\lim_\eb\h$ and
$\z_{\el_2}(e_1)$ are the same. Therefore $\{\ap_i\mid \beta_i\ne 0\}$ are just
the eigenvalues of $h_1$ in $\z_{\el_2}(e_1)_{\Bbb N}$.
Since $[e_1, (\el_2)_i]=(\el_2)_{i{+}1}$ for $i=1,2,\dots$
(see \ref{wond2}(1)), the partition dual to $(\dim(\el_2)_1,\dim(\el_2)_2,
\dots)$ consists of the $h_1$-eigenvalues in $\z_{\el_2}(e_1)_{\Bbb N}$.
But, since the $\Bbb Z$-grading of $\el_2$ is the standard grading associated
with the height of roots
(see the proof of \ref{pr-char}), the dual partition consists also of the
exponents of $\el_2$. This is a classical result of Kostant \cite{k59},
see also \cite[ch.\,4]{CoMc}. \\
This argument is completely symmetric with respect to $e_1$ and $e_2$, because
we do not need the assumption (in the almost principal case) that $e_2$
is Richardson. \\
2. Assume that $\eb$ is almost principal non-integral. Then $\eb$ is principal
in $\g_{{\Bbb Z},{\Bbb Z}}$ (see \cite[2.7]{pCM} or \ref{almeven}(ii))
and we conclude by the first part of the proof.
\end{proof*}%
\subs{Examples} {\sf 1}. $\g={\frak e}_6$. According to \cite[7.6]{wir},
there is a principal nilpotent pair $\eb$ such that $G{\cdot}e_1$ is of type
$\GR{D}{5}$ and $G{\cdot}e_2$ is of type $2\GR{A}{1}$. This means that
$e_1$ (resp. $e_2$) is regular in some Levi subalgebra of type $\GR{D}{5}$
(resp. $2\GR{A}{1}$). We are going to write down explicitly $h_1$ and $h_2$
for this nilpotent pair. Choosing the set of simple roots as in
Theorem~\ref{pr-char},
we see that $h_2$ is dominant and $\el_2$ has to be a standard Levi
subalgebra of type $\GR{D}{5}$. Up to the symmetry of Dynkin diagram,
there is a unique possibility for this:\\
$h_2{=}\Bigl($\begin{tabular}{@{}c@{}}
0\ 0\ \lower3.2ex\vbox{\hbox{0\rule{0ex}{2.4ex}}
\hbox{\hspace{0.4ex}\rule{0ex}{1ex}\rule{0ex}{1ex}}\hbox{0\strut}}\ 0\ 1
\end{tabular}$\Bigr)$. Then
$h_1{=}\Bigl($\begin{tabular}{@{}c@{}}
1\ 1\ \lower3.2ex\vbox{\hbox{1\rule{0ex}{2.4ex}}
\hbox{\hspace{0.4ex}\rule{0ex}{1ex}\rule{0ex}{1ex}}\hbox{1\strut}}\ 1\ $x$
\end{tabular}$\Bigr)$, where $-\cox(\GR{D}{5}){+}1=-7\le x\le 0$.
We identify $h_i$ with the collection $\ap_j(h_i)$, $1\le j\le 6$.
Making use of Ginzburg's results \cite[sect.\,6]{vitya} and some `ad hoc'
arguments, one concludes that $(h_1,h_2)$ can be a characteristic of
a principal pair if and only if $x=-7$. The corresponding
${\Bbb Z}^2$-grading of ${\frak e}_6$ is depicted in Figure 1.
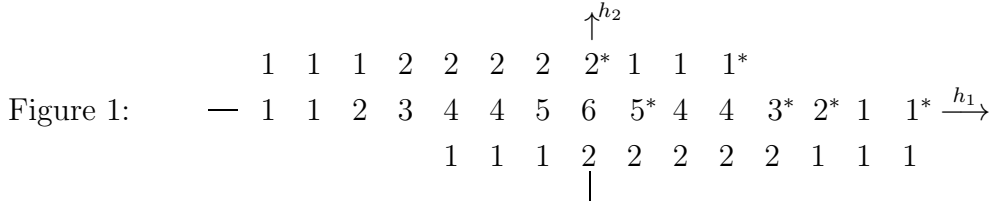
\begin{figure}[ht]
\label{vid6}
\setlength{\unitlength}{0.03in}
\centerline
{
\begin{picture}(130,31)(-65,-13)
\put(-100,0){Figure 1:}
\put(0,0){6}\put(8,0){$5^*$}
\put(16,0){4}\put(24,0){4}\put(32,0){$3^*$}\put(40,0){$2^*$}
\put(48,0){1}\put(56,0){$1^*$}
\put(-8,0){5}
\put(-16,0){4}\put(-24,0){4}\put(-32,0){3}\put(-40,0){2}
\put(-48,0){1}\put(-56,0){1} \put(62,0){$\stackrel{h_1}{\longrightarrow}$}
\put(-65,1.7){\line(1,0){5}}
\put(-8,8){2}
\put(-16,8){2}\put(-24,8){2}\put(-32,8){2}\put(-40,8){1}
\put(-48,8){1}\put(-56,8){1}\put(0,8){$2^*$}\put(8,8){1}
\put(16,8){1}\put(24,8){$1^*$} \put(0,15){$\uparrow$}
\put(2.5,18){\scriptsize $h_2$}
\put(8,-8){2}
\put(16,-8){2}\put(24,-8){2}\put(32,-8){2}\put(40,-8){1}
\put(48,-8){1}\put(56,-8){1}\put(0,-8){2}\put(-8,-8){1}
\put(-16,-8){1}\put(-24,-8){1} \put(1.5,-14){\line(0,1){5}}
\end{picture}
}
\end{figure}
\\[.2ex]
This means, for instance, that $\dim\g_{0,0}=6$,
$\dim\g_{-4,1}=2$ and $\dim\g_{1,0}=5$.
The superscript `*' refers to the eigenspaces containing a one-dimensional
space from $\z_\g(\eb)$.
\\[.6ex]
{\sf 2}. $\g={\frak e}_7$. According to \cite{wir},
there is a principal nilpotent pair $\eb$ such that both $G{\cdot}e_1$ and
$G{\cdot}e_2$ are of type $\GR{A}{4}+\GR{A}{1}$. Here one has several
possibilities for standard Levi subalgebras of type
$\GR{A}{4}+\GR{A}{1}$. Using properties of principal pairs and case-by-case
arguments, one finds that \\
$h_1{=}\Bigl($ \begin{tabular}{@{}c@{}}
$-1$\ 1\ $-4$\ \lower3.4ex\vbox{\hbox{1\rule{0ex}{2.4ex}}
\hbox{\hspace{0.4ex}\rule{0ex}{1ex}\rule{0ex}{1.2ex}}\hbox{1\strut}}\ 1\ 1
\end{tabular}$\Bigr)$ and
$h_2{=}\Bigl($ \begin{tabular}{@{}c@{}}
1\ 0\ 1\ \lower3.1ex\vbox{\hbox{0\rule{0ex}{2.4ex}}
\hbox{\hspace{0.4ex}\rule{0ex}{0.9ex}\rule{0ex}{0.9ex}}\hbox{0\strut}}\ 0\ 0
\end{tabular} $\Bigr)$.
This yields an extremely beautiful ${\Bbb Z}^2$-grading of ${\frak e}_7$,
see Figure 2.
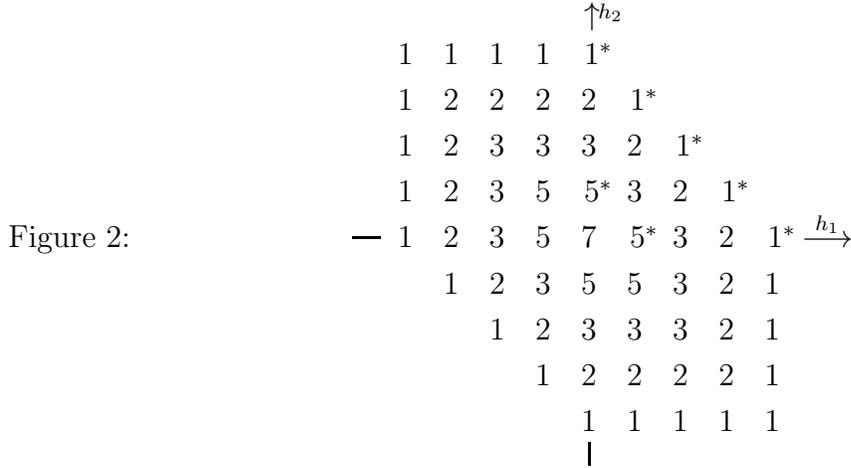
\begin{figure}[ht]
\label{vid7}
\setlength{\unitlength}{0.03in}
\centerline
{
\begin{picture}(130,75)(-65,-36)
\put(-100,0){Figure 2:}
\put(-8,32){1}
\put(-16,32){1}\put(-24,32){1}\put(-32,32){1}
\put(0,32){$1^*$}
 \put(0,39){$\uparrow$} \put(2.5,40){\scriptsize $h_2$}
\put(-8,24){2}
\put(-16,24){2}\put(-24,24){2}\put(-32,24){1}
\put(0,24){2}\put(8,24){$1^*$}
\put(-8,16){3}
\put(-16,16){3}\put(-24,16){2}\put(-32,16){1}
\put(0,16){3}\put(8,16){2}
\put(16,16){$1^*$}
\put(-8,8){5}
\put(-16,8){3}\put(-24,8){2}\put(-32,8){1}
\put(0,8){$5^*$}\put(8,8){3}
\put(16,8){2}\put(24,8){$1^*$}
\put(0,0){7}\put(8,0){$5^*$}
\put(16,0){3}\put(24,0){2}\put(32,0){$1^*$}
\put(-8,0){5}
\put(-16,0){3}\put(-24,0){2}\put(-32,0){1}
\put(38,0){$\stackrel{h_1}{\longrightarrow}$}
\put(-40,1.7){\line(1,0){5}}
\put(8,-8){5}
\put(16,-8){3}\put(24,-8){2}\put(32,-8){1}
\put(0,-8){5}\put(-8,-8){3}
\put(-16,-8){2}\put(-24,-8){1}
\put(8,-16){3}
\put(16,-16){3}\put(24,-16){2}\put(32,-16){1}
\put(0,-16){3}\put(-8,-16){2}\put(-16,-16){1}
\put(8,-24){2}
\put(16,-24){2}\put(24,-24){2}\put(32,-24){1}
\put(0,-24){2}\put(-8,-24){1}
\put(8,-32){1}
\put(16,-32){1}\put(24,-32){1}\put(32,-32){1}
\put(0,-32){1}
\put(1.4,-38){\line(0,1){4}}
\end{picture}
}
\end{figure}
\\ 
{\sf 3.} In \cite[2.15]{pCM}, we described a series of integral almost
principal pairs in ${\frak sp}_{4n}$. The orbit $G{\cdot}e_1$ corresponds
to the partition $(2n,2n)$ and is Richardson whereas $G{\cdot}e_2$ corresponds
to $(2^{2n-1},1,1)$ and is not Richardson. The formulas for $\hb$ imply that
$\el_2=\GR{A}{2n-1}$ and $\el_1=n\GR{A}{1}$. These formulas are such that
neither of $h_i$ is dominant with respect to the
standard set of simple roots.
Choosing $\Pi$ adapted to $h_2$, as in Theorem \ref{pr-char}, we obtain
$h_1=(1\ 1\ldots 1\ -2n)$, $h_2=(0\ 0\ldots 0\ 1)$.
Since $-2n {<} -2n{+}1=-\cox(\GR{A}{2n-1}){+}1$, we see that \ref{pr-char}(iii)
does not hold here. But choosing $\Pi$ adapted to $h_1$, we obtain
$h_1=(1\ 0\ldots 1\ 0)$, $h_2=(-1\ 1\ldots -1\ 1)$, which correlates with
the fact that $\cox(\GR{A}{1})=2$.
This shows the condition of being Richardson for $G{\cdot}e_2$ is essential
in \ref{pr-char}(iii).

\subs{Classification problems\nopagebreak}%
A usual intention of a mathematician is to classify something,
especially those objects whose number is finite. Turning back, we see
at least four groups of objects that could be classified: \\
\indent $\bullet$ the $G$-orbits of characteristics of nilpotent pairs; \\
\indent $\bullet$ the $G$-orbits of wonderful nilpotent pairs; \\
\indent $\bullet$ the $G$-orbits of even and almost
even nilpotent pairs; \\
\indent $\bullet$ the $G$-orbits of principal and almost
principal nilpotent pairs. \\
The last group is the smallest one, and a complete description for it will
be given in \cite{wir2}. But obtaining a classification of all characteristics
requires at least a better understanding of possible fractional eigenspaces
and, at the moment, I have no idea for this.

\vno{3} \indent
{\footnotesize
\parbox{255pt}{%
{\it Math. Department \\
MIREA \\
prosp. Vernadskogo, 78 \\
Moscow 117454 \quad Russia} \\ dmitri@panyushev.mccme.ru \qquad
panyush@dpa.msk.ru }
}

\begin{thebibliography}{Pa95a}
{\footnotesize

\bibitem[Br89]{ranee} {\sc Brylinski, R.K.}: Limits of weight spaces, Lusztig's
$q$-analogs, and fiberings of adjoint orbits, {\it J. Amer. Math. Soc.}
{\bf 2}(1989), 517--533.

\bibitem[CM93]{CoMc}
{\sc Collingwood, D.H.; McGovern, W.M.}: "Nilpotent orbits in semisimple
  Lie algebras", New York: Van Nostrand Reinhold, 1993.

\bibitem[Dy52a]{EBD}
{\sc Dynkin, E.B.}: Semisimple subalgebras of semisimple Lie algebras,
{\it Matem. Sbornik} {\bf 30}(1952), {\rus N0}\,2, 349--462
(Russian). English translation: {\it Amer. Math. Soc. Transl.} II~Ser.,
{\bf 6}~(1957), 111--244.

\bibitem[Dy52b]{max}
{\sc Dynkin, E.B.}: The maximal subgroups of the classical groups,
{\it Trudy Mosk. Matem. Obsch.} {\bf 1}(1952), 39--166 (Russian).
English translation: {\it Amer. Math. Soc. Transl.} II~Ser.,
{\bf 6}~(1957), 245--378.

\bibitem[El75]{El75}
{\sc Elashvili, A.G.}: The centralizers of nilpotent elements in semisimple Lie
algebras, {\it Trudy Tbiliss. Matem. Inst. Akad. Nauk Gruzin. SSR}
{\bf 46}(1975), 109--132 (Russian).

\bibitem[EP99]{wir}
{\sc Elashvili, A.G.; Panyushev, D.}:
Towards a classification of principal nilpotent pairs,
{\it Appendix} to [Gi99].

\bibitem[ElPa]{wir2}
{\sc Elashvili, A.G.; Panyushev, D.}:
A classification of principal nilpotent pairs and related problems,
{\it in preparation}.

\bibitem[Gi99]{vitya} {\sc Ginzburg, V.}: Principal nilpotent pairs
in a semisimple Lie algebra I, {\it Preprint}\\ math.AG/9903059.

\bibitem[Ko59]{k59} {\sc Kostant, B.}: The principal three-dimensional
subgroup and the Betti numbers of a complex simple Lie group,
{\it Amer. J. Math.} {\bf 81}(1959), 973--1032.

\bibitem[Pa99]{pCM} {\sc Panyushev, D.}: Nilpotent pairs,
dual pairs, and sheets, {\it Preprint\/} math.AG/9904014.

\bibitem[Conj]{SS}
{\sc Springer, T.A.; Steinberg, R.}: Conjugacy classes, In:
{\it "Seminar on algebraic groups and related finite groups".} Lecture notes
in Math. {\bf 131}, pp.167--266,  Berlin
Heidelberg New York: Springer 1970.

\bibitem[Vi79]{vi79}
{\sc Vinberg, E.B.}: Classification of homogeneous nilpotent elements
of a semisimple graded Lie algebra, In: {\it "Trudy seminara po vect.
i tenz. analizu"}, vol.\,19, pp.~155--177. Moscow: MGU 1979 (Russian).
English translation:
{\it Selecta Math. Sovietica} {\bf 6}(1987), 15--35.

\bibitem[VGO90]{t41}
{\sc Vinberg, E.B.; Gorbatsevich, V.V.; Onishchik, A.L.}:
``Structure of Lie groups and Lie algerbas", In: {\it
Sovremennye problemy matematiki. Fundamentalnye napravleniya}, t.\,41.
Moskwa: VINITI 1990 (Russian).
English translation in: V.V.\,Gorbatsevich, A.L.\,Onishchik, E.B.\,Vinberg,
Lie Groups and Lie Algebras III
(Encyclopaedia Math. Sci., vol.~41) Berlin Heidelberg New York:
Springer 1994.

\bibitem[VO88]{vion}
{\sc Vinberg, E.B.; Onishchik, A.L.}: ``Seminar on Lie groups and
algebraic groups",  Moscow: Nauka 1988 (Russian).
English translation: {\sc Onishchik, A.L.; Vinberg, E.B.}:
Lie groups and algebraic groups. Berlin Heidelberg
New York: Springer 1990)
}
\end{thebibliography}
\end{document}